\documentclass[11pt]{amsart} 

\usepackage{amssymb} 
\usepackage{enumerate} 
\usepackage{paralist}
\usepackage{hyperref} 
\usepackage{tikz} 
\usepackage{amsmath}

\theoremstyle{definition} 
\newtheorem{theorem}{Theorem}[section]
\newtheorem{lemma}[theorem]{Lemma}

\theoremstyle{definition} 
\newtheorem{definition}[theorem]{Definition}

\title{Linear orders and the real line}
\author{Aksel Ozer, Trey Smith}

\begin{document}


\maketitle

\section{Abstract}
\bigskip
We define both dense and separable linear orders. We state and prove unique order-theoretic characterizations of both ($\mathbb{R},<$) and ($\mathbb{Q}, <$) due to Cantor. We describe Suslin's Problem and give several equivalent statements. We then discuss Aronszajn Trees and their relation to Suslin Trees. Finally we briefly touch on the independence of Suslin's Problem from ZFC and the possible existence of Suslin Lines/Trees.

\bigskip

\section{Background}
\bigskip
In 1897, Georg Cantor characterized the order-types of the rational and real numbers as a part of his theory of transfinite sets. In 1920 a list of problems was published in the first volume of the Polish journal \emph{Fundamenta Mathematicae}. The third problem on the list was posthumously attributed to Mikhail Suslin \cite{kanamori_2011}. This problem asked whether one of the conditions in Cantor's characterization of the real numbers could be weakened \cite{suslin_1920}. This relatively benign question has come to influence decades of set theory. This list of problems had to do with possible consequences of the Continuum Hypothesis and issues in the newly formed descriptive set theory. While the other problems on the list were solved, Suslin's problem remained without a solution and began to grow in significance as a result. Suslin himself continued to contribute to the theory of analytic sets (sometimes called Suslin sets). Unfortunately, in 1919 he met an untimely end due to typhus in Moscow, following the Russian Civil War. One of most significant results about Suslin's Problem, one that would inspire much future research, was the result of Ronald Jensen. Which stated that for G\"odel’s constructible universe L, if V = L, then there is a Suslin tree, which we will see solves Suslin's Problem. Jensen famously found a combinatorial principle that  implies that there is a Suslin tree \cite{JENSEN_1972}. This would lead to generalizations to higher cardinals, and also begin Jensen’s study of pivotal combinatorial principles holding in L, principles that would establish many of the propositions of combinatorial set theory. Another significant result was that of Solovay and Tennenbaum which built on Cohen's method of forcing, and was later generalized by Donald Martin and Solovay into Martin's Axiom which implies the Suslin Hypothesis \cite{martin_solovay_1970}.

\bigskip

\newpage

\section{Cantor's Theorem}

\begin{definition}
    The theory of linear orders

The properties that characterize linear orders are listed below:\\
\begin{enumerate}
\item(Axiom 1) Linear order: $\forall x \forall y(x<y \vee y<x \vee y=x)$
for any two elements on the line, they are either equal or one is greater than the other\\
\item(Axiom 2) Linear order: $\forall x\neg( x<x)$
 No number can be less than itself, also known as the property of irreflexivity\\ 
\item(Axiom 3) Transitivity: $\forall x \forall y \forall z(x<y \wedge y<z \rightarrow x<z)$\\

\end{enumerate}
\end{definition}

Examples of dense linear orders are $(\mathbb{Q},<)$ and $(\mathbb{R},<)$.\\

\begin{definition}[Dense Linear Order or "DLO"]

  A linear order $(P, <)$ is dense if for all $a<b$ there exists a $c$ such that $a<c<b$.
  \label{def:Dense Linear Order} 
\end{definition}

\bigskip

We now give a proof of Cantor's unique characterization of the rational numbers \cite{jech_1997}.

\bigskip

\begin{theorem}[Cantor’s Back-and-Forth Theorem]
$(\,\mathbb{Q}, <)\,$ is the unique dense linear order without endpoints.\\\\
\emph{Proof.}  Let $A=\left\{a_{n}: n \in \boldsymbol{N}\right\}$ and let $B=\left\{b_{n}: n \in \boldsymbol{N}\right\}$ be two
unbounded, dense, linear orders.\\ We construct an isomorphism $f: A \rightarrow B$ in the following way: We first define partial isomorphisms $f\left(a_{0}\right)$, then $f^{-1}\left(b_{0}\right)$, then $f\left(a_{1}\right)$, then $f^{-1}\left(b_{1}\right)$, etc., so as to keep $f$ order-preserving. For example, to define $f\left(a_{n}\right)$, if it is not yet defined, we let $f\left(a_{n}\right)=b_{m}$ where $m$ is the smallest index such that $f$ remains order-preserving (This will $m$ always exists because $f$ has been defined for only finitely many $a \in A$, and because $B$ is dense and unbounded). So $f$ is an order-preserving bijection. Thus $A\simeq B$ \hfill\(\Box\)
\end{theorem}
Thus any nonempty, countable, unbounded, dense linear order is isomorphic to the rational numbers.\\

\begin{definition}[Separable Linear Order]
    A linear order $(P, <)$ is separable if there exists a countable set C $\subseteq$ P such that for all $a, b \in P$ with $a < b$, there is a $c \in C$ such that $a < c < b$, i.e. C is dense in P.\\
\end{definition}

\begin{definition}[Complete Linear Order]
    A linear order $(P, <)$ is said to be complete if every nonempty subset that has an upper bound, has a least upper bound.\\
\end{definition}

\begin{definition}[Unbounded]
    A linear order $(P, <)$ is unbounded if for all $a \in P$ there are $b, c \in P$ such that $b < a < c$.
\end{definition}

We now give a proof of Cantor's unique characterization of the real numbers \cite{jech_1997}.

\bigskip

\begin{theorem}[Cantor's characterization of the Real Numbers] ($\mathbb{R},<$) is the unique complete linear ordering that has a countable dense subset isomorphic to ($\mathbb{Q}, <$).\\\\
\emph{Proof.} Let $D$ and $D'$ be two complete dense
unbounded linearly ordered sets, let $P$ and $P'$ be dense in $D$ and $D'$, respectively. By Theorem 3.3 we can give an isomorphism $I$ of $P$ onto $P'$. Then $I$ can be uniquely extended to an isomorphism $I^\ast$ of $D$ and $D'$: For $x \in D$, let $I^\ast(x) = \sup\{I(p) : p \in P$ and $p \leq x\}$. $I^\ast$ is clearly injective since if two elements share the same preceding elements, then they are the same element. $I^\ast$ is also surjective because of the density of $P'$  in $D'$. Finally since $I$ is order preserving, $I^\ast$ will be as well. Thus $D\simeq D'$. \hfill\(\Box\)

\end{theorem}

\bigskip

Thus any complete, dense, unbounded, separable linear order is isomorphic to the real line. Now that we have Cantor's unique characterizations of the rationals and the reals, we turn our attention to Suslin's Problem.

\bigskip

\section{Suslin's Problem}
\begin{definition}[The countable chain condition]
    A linear order $(P, <)$ satisfies the countable chain condition if every family of pairwise disjoint, non-empty open intervals
$(a, b) = \{c \in P : a < c < b\}$ is countable.

\end{definition}
\begin{figure}
\begin{center}
\includegraphics[scale=1.5]{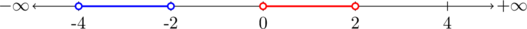}
    \caption{Disjoint open intervals on the real line}
\end{center} 
\end{figure}
\bigskip

The Real Line is dense, unbounded, complete, separable, and satisfies the countable chain condition since every collection of mutually disjoint, nonempty open intervals in $\mathbb{R}$ is countable. \textbf{Suslin's Problem} asks whether separable can be replaced by the countable chain condition. In other words, is every linear order which is dense, complete, unbounded, and satisfies the countable chain condition, isomorphic to $\mathbb{R}$ ? \textbf{Suslin's Hypothesis} is the affirmative answer to Suslin's Problem \cite{kanamori_2011}.\\\\

\begin{definition}[Suslin Lines]
A Suslin Line is a dense, unbounded linear order that satisfies the countable chain condition but is not separable.
\end{definition}
Thus Suslin's Problem is equivalent to the problem of whether Suslin Lines exist, since if they do exist, Suslin's Hypothesis will be true. Now we will discuss Suslin Trees and their equivalence to Suslin Lines.

\bigskip

\begin{definition}[Tree]
    A tree is a partially ordered set $(T,<)$ with the property
that for each $x \in T$ , the set $\{y : y<x\}$ of all predecessors of $x$ is well-ordered by $<$.
\end{definition}
    \[o(x) = \mbox {the order-type of } \{y : y<x\},\]
    \[\alpha \mbox{th level} = \{x : o(x) = \alpha\},\]
    \[\mbox{height}(T) = \sup\{o(x)+1: x \in T \}.\]
\vspace{4mm} 
    
    A \emph{branch} in $T$ is a maximal linearly ordered subset of $T$. The length of a branch $b$ is the order-type of $b$. An $\alpha$-branch is a branch of length $\alpha$.\\

\begin{figure}  
\begin{center}
    \includegraphics[scale=1.7]{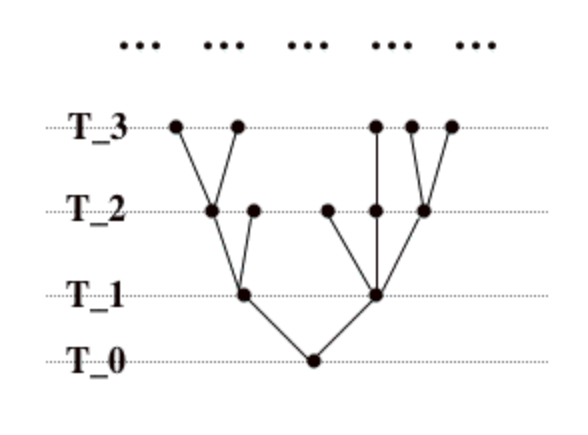}
        \caption{$\mbox{A tree with levels } T_0, T_1, T_2, \mbox{and }T_3$}
\end{center}
\end{figure}
\bigskip    
    Let T be a tree. An \emph{antichain} in $T$ is a set $A \subset T$ such that any two distinct elements $x, y$ of $A$ are incomparable, i.e., neither $x<y$ nor $y<x$.

\bigskip

\begin{definition}[The first uncountable ordinal]
$\omega_1$ is the first uncountable ordinal, which has cardinality $|\mathbb{Q}|<|\omega_1|\leq|\mathbb{R}|$. 
\end{definition}

This can be seen by showing that every countable ordinal embeds into ($\mathbb{Q}, <$) and thus into ($\mathbb{R},<$).\\

\begin{theorem}[Every countable linear order embeds into ($\mathbb{Q}, <$)] There is an isomorphism between every countable linear order and a subset of $\mathbb{Q}$\\\\
\emph{Proof.} Since $A$ is countable, we can arrange it in a sequence $a_{1}, a_{2}, a_{3}, \ldots$ We describe a
procedure to define $f\left(a_{i}\right)$ for each $a_{i}$ in turn. Let $f\left(a_{1}\right)$ be any rational. Suppose we have
defined $f\left(a_{1}\right), f\left(a_{2}\right), \ldots, f\left(a_{n}\right)$ in such a way that all order relations are preserved (that is, for
all $i, j \leq n, a_{i} \leq a_{j}$ if and only if $\left.f\left(a_{i}\right) \leq f\left(a_{j}\right)\right) .$ We want to define $f$ on $a_{n+1} \in A$. Partition the set $\left\{a_{1}, \ldots, a_{n}\right\}$ into two subsets:
$$
X=\left\{a_{i}: i \leq n \text { and } a_{i}<a_{n+1}\right\}, Y=\left\{a_{i}: i \leq n \text { and } a_{i}>a_{n+1}\right\}
$$
In $\mathbb{Q}$, every element of $f(X)$ is smaller than every element of $f(Y) .$ Choose $q$ strictly larger than the elements of $f(X)$ and strictly smaller than the elements of $f(Y)$. For each $i \leq n$, the
relationship between $q$ and $f\left(a_{i}\right)$ is the same as the relationship between $a_{n+1}$ and $a_{i} .$

Therefore, letting $f\left(a_{n+1}\right)=q$, we have extended the function to one more element in such a
way that all order relations are preserved. The resulting function defined on all of $A$ is thus an
isomorphism from $A$ to the range of $f$.\hfill\(\Box\)
    
\end{theorem}

\bigskip

\begin{definition}[Suslin Trees]
    A tree $T$ is a Suslin Tree if:
\begin{enumerate}
    \item the height of $T$ is $\omega_1$;
    \item every branch in $T$ is at most countable;
    \item every antichain in $T$ is at most countable;
\end{enumerate}
\end{definition}
    In order to formulate Suslin’s Problem in terms of trees it is helpful to think about Suslin Trees that are called \emph{normal} \cite{jech_1997}. Let $\alpha$ be an ordinal number, $\alpha \leq \omega_1$. A normal $\alpha$-tree is a tree $T$ with the following properties:\\
\begin{enumerate}
    \item height($T$) = $\alpha$;
    \item $T$ has a unique least point (the "root" of the tree, i.e. $T_0$ in Figure 2);
    \item each level of $T$ is at most countable;
    \item if $x$ is not maximal in $T$ , then there are infinitely many $y>x$ at the next level (immediate successors of $x$);
    \item for each $x \in T$ there is some $y>x$ at each higher level less than $\alpha$;
    \item if $\beta<\alpha$ is a limit ordinal and $x, y$ are both at level $\beta$ and if $\{z : z<x\} = \{z : z<y\}, \mbox{then } x = y.$
\end{enumerate}

\bigskip

\begin{lemma}\emph{If there exists a Suslin Tree then there exists a normal Suslin Tree.}\\\\
    \emph{Proof. }Let $T$ be a Suslin tree. $T$ has height $\omega_1$, and each level of $T$ is countable. We first get rid of all points $x \in T$ where $T_x = \{y \in T : y \geq x\}$ is at most countable, and let $T_1 = \{x \in T : T_x$ is uncountable$\}$. We can see that if $x \in T_1$ and $\alpha > o(x)$, then $|T_y| = \aleph_1$ for some $y>x$ at level $\alpha$. Hence $T_1$ satisfies condition (5). Next, we satisfy property (6): For every chain $C = \{z : z<y\}$ in $T_1$ of limit length we add an extra node $a_C$ and stipulate that $z<a_C$ for all $z \in C$, and $a_C < x$ for every $x$ such that $x>z$ for all $z \in C$.\\\\ 
    The resulting tree $T_2$ satisfies (3), (5) and (6). For each $x \in T_2$ there are uncountably many branching points $z>x$, i.e., points that have at least two immediate successors (because there is no uncountable chain and $T_2$ satisfies (5)). The tree $T_3 = \{x \in T_2 : x \mbox{ is a branching point}\}$ satisfies (3), (5) and (6) and each $x \in T_3$ is a branching point. To get property (4), let $T_4$ consists of all $z \in T_3$ at limit levels of $T_3$. The tree $T_4$ satisfies (1), (3), (4), and (5); and then $T_5 \subset T_4$ satisfying (2) as well is easily obtained. \hfill\(\Box\)
\end{lemma}

\begin{lemma}\emph{There exists a Suslin Line iff there exists a Suslin Tree}\\\\
    \emph{Proof. }($\Rightarrow$) Let $S$ be a Suslin line. We will construct a Suslin tree. The tree will be made of nondegenerate, closed intervals on $S$. The partial ordering of $T$ is by inverse inclusion: If $I,J \in T$ , then $I \leq J$ iff $I \supset J$. The collection $T$ of intervals is constructed by induction on $\alpha<\omega_1$. We let $I_0 = [a_0, b_0]$ be arbitrary (such that $a_0 < b_0$). Having constructed $I_\beta, \beta<\alpha$, we consider the countable set $C = \{a_\beta : \beta<\alpha\}\cup\{b_\beta : \beta<\alpha\}$ of endpoints of the intervals $I_\beta, \beta<\alpha$. Since $S$ is a Suslin line, $C$ is not dense in $S$ and so there exists an interval $[a, b]$ disjoint from $C$; we pick some such $[a_\alpha, b_\alpha] = I_\alpha$. The set $T = \{I_\alpha : \alpha<\omega_1\}$ is uncountable and partially ordered by $\supset$. If $\alpha<\beta$, then either $I_\alpha \supset I_\beta$ or $I_\alpha$ is disjoint from $I_\beta$. It follows that for each $\alpha, \{I \in T : I \supset I_\alpha\}$ is well-ordered by $\supset$ and thus $T$ is a tree.\\\\ 
    We now show that $T$ has no uncountable branches and no uncountable antichains. Then it is clear that the height of $T$ is at most $\omega_1$ and since every level is an antichain and $T$ is uncountable, we have height($T$) $= \omega_1$. If $I,J \in T$ are incomparable, then they are disjoint intervals of $S$ and since $S$ satisfies the countable chain condition, every antichain in $T$ is at most countable. To show that $T$ has no uncountable branch, we note first that if $b$ is a branch of length $\omega_1$, then the left endpoints of the intervals $I \in B$ form an increasing sequence $\{x_\alpha : \alpha<\omega_1\}$ of points of $S$. It is clear that the intervals $(\,x_\alpha, x_{\alpha+1})\,$ for $\alpha<\omega_1$, form a disjoint uncountable collection of open intervals in $S$, contrary to the assumption that S satisfies the countable chain condition.\\\\
    
    ($\Leftarrow$) Let $T$ be a normal Suslin tree. The line $S$ will consist of all branches of $T$ (which are all countable). Each $x \in T$ has countably many immediate successors, and we order these successors as rational numbers. Then we order the elements of $S$ lexicographically: If $\alpha$ is the least level where two branches $a, b \in S$ differ, then $\alpha$ is a successor ordinal and the points $a_\alpha \in a$ and $b_\alpha \in b$ are both successors of the same point at level $\alpha-1$; we let $a<b$ or $b<a$ according to whether $a_\alpha < b_\alpha$ or $b_\alpha < a_\alpha$. It is easy to see that S is linearly ordered and dense. If $(\,a, b)\,$ is an open interval in $S$, then one can find $x \in T$ such that $I_x \subset (\,a, b)\,$, where $I_x$ is the interval $I_x = \{c \in S : x \in c\}$. And if $I_x$ and $I_y$ are disjoint, then $x$ and $y$ are incomparable points of $T$ .\\\\ 
    Thus every disjoint collection of open intervals of $S$ must be at most countable, and so $S$ satisfies the countable chain condition. The line $S$ is not separable: If $C$ is a countable set of branches of $T$ , let $\alpha$ be a countable ordinal bigger than the length of any branches $b \in C$. Then if $x$ is any point at level $\alpha$, the interval $I_x$ does not contain any $b \in C$, and so $C$ is not dense in $S$.\hfill\(\Box\)
    \end{lemma}

\bigskip

Thus the existence of a Suslin Tree is equivalent to Suslin's Problem since we have seen that the existence of a Suslin Line is equivalent to Suslin's Hypothesis. We can almost give a counterexample to Suslin's Problem in ZFC called an Aronszajn Tree.\\
\begin{definition}[Aronszajn Tree]
    An Aronszajn tree is a tree of height $\omega_1$ all of whose levels are at most countable and which has no uncountable branches. 
\end{definition}

Thus every Suslin Tree is an Aronszajn Tree, however not all Aronszjan Trees are Suslin Trees. We will now give an example of an Aronszajn Tree which is not a Suslin Tree.\\

\begin{theorem}[Aronszajn Tree]
    There exists an Aronszajn Tree.\\\\
    \emph{Proof.} We construct a tree T whose elements are some bounded increasing transfinite sequences of rational numbers. If $x, y \in T$ are two such sequences, then we let $x \leq y$ just in case y extends $x$, i.e., $x \subset y$. Also, if $y \in T$ and $x$ is an initial segment of $y$, then we let $x \in T$ ; thus the $\alpha$th level of $T$ will consist of all those $x \in T$ whose length is $\alpha$. It is clear that an uncountable branch in $T$ would yield an increasing $\omega_1$-sequence of rational numbers, which is impossible. Thus $T$ will be an Aronszajn tree, provided we arrange that $T$ has $\aleph_1$ levels, all of them at most countable. We construct $T$ by induction on levels. For each $\alpha<\omega_1$ we construct a set $U_\alpha$ of increasing $\alpha$-sequences of rationals; $U_\alpha$ will be the $\alpha$th level of $T$ . We construct the $U_\alpha$ so that for each $\alpha$, $|U_\alpha|\leq\aleph_0$, and that: 
\bigskip    
\begin{enumerate}
    
    \item$\mbox{ For each } \beta<\alpha, \mbox{ each } x \in U_\beta \mbox{ and each } q \geq \sup x \mbox{ there is }y \in U_\alpha\\ \mbox{ such}\mbox{ that }x \subset y \mbox{ and }q \geq \sup y.$

\end{enumerate}
\bigskip    
    This condition enables us to continue the construction at limit steps. To start, we let $U_0 = \{\emptyset\}$. The successor steps of the construction are also fairly easy. Given $U_\alpha$, we let $U_{\alpha+1}$ be the set of all extensions $x\frown r$ of sequences in $U_\alpha$ such that $r > \sup x$. It is clear that since $U_\alpha$ satisfies condition (1), $U{\alpha+1}$ satisfies it also (for $\alpha$ + 1), and it is equally clear that $U_{\alpha+1}$ is at most countable. Thus let $\alpha$ be a limit ordinal ($\alpha<\omega_1$) and assume that we have constructed all levels $U_\gamma, \gamma<\alpha$, of $T$ below $\alpha$, and that all the $U_\gamma$ satisfy (1); we shall construct $U_\alpha$. The points $x \in T$ below level $\alpha$ form a (normal) tree $T_\alpha$ of length $\alpha$. We claim that $T_\alpha$ has the following property:
\bigskip
\begin{enumerate}
    \setcounter{enumi}{1}
    \item For each $x \in T_\alpha$ and each $q > \sup x$ there is an increasing $\alpha$-sequence of rationals $y$ such that $x \subset y$ and $q \geq \sup y$ and that $y\restriction\beta \in T_\alpha$ for all $\beta<\alpha$.

\end{enumerate}
\bigskip
    The last condition means that $\{y\restriction\beta : \beta<\alpha\}$ is a branch in $T_\alpha$. To prove the claim, we let $\alpha_n$, $n = 0, 1, . . .$ , be an increasing sequence of ordinals such that $x \in U_{\alpha_0}$ and $\lim_n \alpha_n = \alpha$, and let $\{q_n\}_{n=0}^{\infty}$ be an increasing sequence of rational numbers such that $q_0 > \sup x$ and $\lim_n q_n \leq q$. Repeatedly using condition (1), for all $\alpha_n $ $(\,n = 0, 1, . . . )\,$, we can construct a sequence $y_0 \subset y_1 \subset ... \subset y_n ... $ such that $y_0 = x, y_n \in U_{\alpha_n}$ , and $\sup y_n \leq q_n$ for each $n$. Then we let $y = U_{n=0}^{\infty} y_n$; clearly, $y$ satisfies condition (2). Now we construct $U_\alpha$ as follows: For each $x \in T_\alpha$ and each rational number $q$ such that $q > \sup x$, we choose a branch $y$ in $T_\alpha$ that satisfies (2), and let $U_\alpha$ consists of all these $y : \alpha \rightarrow Q$. The set $U_\alpha$ so constructed is countable and satisfies condition (1). Then $T = \cup_{\alpha<\omega_1} U_\alpha$ is an Aronszajn tree. \hfill\(\Box\)
\end{theorem}

This is a particular kind of Aronszajn Tree known as a \emph{Special Aronszajn Tree}, meaning that we can give a strictly increasing function from the tree to the rational numbers. This also means it is not a Suslin Tree since it will have an uncountable antichain.

\bigskip

\section{Independence from ZFC}
Suslin's Problem has been shown to be independent from ZFC, the usual Zermelo-Fraenkel axioms of set theory with the Axiom of Choice \cite{kunen_2006}. This means that there are models of ZFC where the Suslin Hypothesis is true, and models where the Suslin Hypothesis is false. Interestingly, the proofs for both kinds of models involve the method of forcing. Forcing involves starting with a given model M called the ground model, and extending it by adjoining an new object G. The new model M[G] is basically a minimal possible collection of sets that includes M, contains G, and also satisfies ZFC \cite{jech_1997(2)}. In 1967, Thomas Jech used forcing methods to give a model of ZFC where Suslin Lines exist \cite{Jech_1967}. This was also proved independently by Stanley Tennenbaum likewise using the method of forcing \cite{Tennenbaum_60}. In 1971, Tennenbaum and a mathematician named Robert Solovay used forcing to give a model of ZFC where Suslin Lines do not exist \cite{Solovay1971IteratedCE}. Thus, like the Continuum Hypothesis and $\aleph_1$, the existence of Suslin Lines is an open question. It has been shown by Jensen that the Axiom of Constructibility implies the existence of Suslin Lines \cite{JENSEN_1972}. However, this axiom is not yet widely accepted among set theorists. As mentioned in the background, Solovay and Martin generalized these forcing techniques to propose Martins Axiom, which also implies Suslin's Hypothesis \cite{kunen_2006}. Martin's Axiom roughly states that cardinals smaller than $|\mathbb{R}| = \mathfrak{C}$ behave like $\aleph_0$ in the sense that given a partially ordered set ($P, <$) with cardinality $\lambda < \mathfrak{C}$ we can construct a specific kind of filter (a special subset) on $P$ \cite{martin_solovay_1970}. Unfortunately, Martin's Axiom is only interesting if the Continuum Hypothesis is false, otherwise there are no cardinals between $\aleph_0$ and $\mathfrak{C}$.

\newpage

\bibliographystyle{plain}
\bibliography{Citations}

\end{document}